\input amstex\documentstyle{amsppt}  
\pagewidth{12.5cm}\pageheight{19cm}\magnification\magstep1  
\topmatter
\title On the Satake isomorphism\endtitle
\author G. Lusztig\endauthor
\address{Department of Mathematics, M.I.T., Cambridge, MA 02139; Institute for Advanced Study,
Princeton, NJ 08450}\endaddress
\thanks{Supported by NSF grant DMS-1855773.}\endthanks
\endtopmatter   
\document

\define\pos{\text{\rm pos}}

\define\mpb{\medpagebreak}

\define\frl{\forall}

\define\si{\sim}

\define\sqc{\sqcup}

\define\qua{\quad}

\define\op{\oplus}
   
\define\part{\partial}

\define\n{\notin}

\define\m{\mapsto}
\define\do{\dots}

\define\sub{\subset}    

\define\T{\times}
\define\ti{\tilde}
\define\nl{\newline}
\redefine\i{^{-1}}
\define\fra{\frac}
\define\un{\underline}

\define\ot{\otimes}
\define\bbq{\bar{\QQ}_l}

\define\tr{\text{\rm tr}}

\define\a{\alpha}
\redefine\b{\beta}

\redefine\d{\delta}

\redefine\o{\omega}
\define\p{\pi}

\define\r{\rho}
\define\s{\sigma}
\redefine\t{\tau}
\define\th{\theta}

\define\x{\xi}

\define\boc{\bold c}

\define\kk{\bold k}

\redefine\AA{\bold A}

\define\CC{\bold C}

\define\FF{\bold F}

\define\HH{\bold H}

\define\LL{\bold L}

\define\NN{\bold N}

\define\QQ{\bold Q}
\define\RR{\bold R}
\define\SS{\bold S}

\define\VV{\bold V}

\define\ZZ{\bold Z}

\define\ca{\Cal A}
\define\cb{\Cal B}

\define\cf{\Cal F}

\define\ch{\Cal H}

\define\co{\Cal O}
\define\cp{\Cal P}

\define\cv{\Cal V}
\define\cw{\Cal W}
\define\cz{\Cal Z}

\define\fg{\frak g}

\define\tW{\ti W}

\define\sha{\sharp}

\define\bP{\bar P}

\define\che{\check}
\define\cha{\che{\a}}

\head Introduction\endhead

\subhead 0.1\endsubhead
Let $H_q$ be the affine Hecke algebra over $\CC$
(with equal parameters $q$, a prime power) associated to
an affine Weyl group $W$ (defined in terms of the dual $G^*$ of an adjoint group $G$).
Let $H_{0,q}$ be the Hecke algebra over $\CC$ (with equal parameters $q$)
associated to the corresponding finite Weyl group $W_0\sub W$.
Let $H_q^{sph}$ the vector subspace of $H_q$
consisting of elements which are eigenvectors for the left and right multiplication by $H_{0,q}$
with eigenvalue defined by the one dimensional representation of $H_{0,q}$ corresponding to the
unit representation of $W_0$. Then $H_q^{sph}$ is an algebra for the product
$f*f'=(\sum_{w\in W_0}q^{|w|})\i ff'$
where $ff'$ is the product in $H_q$ and $||$ is the standard length function on $W_0$.
Let $Q$ be the group of translations in $W$.
The classical Satake isomorphism states that the algebra $H_q^{sph}$ is isomorphic
to the algebra of $W_0$-invariants in the group algebra $\CC[Q]$.
In \cite{L83} we gave a refinement of this isomorphism in which
the basis of $\CC[Q]$ formed by the irreducible representations
of a semisimple group with Weyl group $W_0$ and for which $Q$ is the lattice of roots
corresponds to a basis $\b$ of $H_q^{sph}$
formed by certain elements of the basis \cite{KL79} of $H_q$, suitably normalized.
This shows in particular that

(a) the structure constants of the algebra $H_q^{sph}$ with respect to $\b$ are integers
independent of $q$.
\nl
This is the starting point of ``geometric Satake
equivalence'' (which we do not discuss in this paper).

\subhead 0.2\endsubhead
In this paper we show (see 1.5)
that the structure constants in 0.1(a) can be interpreted as structure
constants for a certain subring $J_*$
of the $J$-ring attached to $W$ with respect to the standard basis of $J_*$.
(We actually prove a more general statement involving a weight function on $W$.)
This gives a new (and simpler) proof of 0.1(a). We also give another approach to 0.1(a)
based on the character formula for simple rational modules of a semisimple group in
characteristic $p\gg0$.
 At the time when \cite{L83} was written, this character formula
was only conjectured and providing evidence for the conjecture was one of the
motivations which led the author to \cite{L83}.
We also state an extension of that character formula to
disconnected groups.

\subhead 0.3\endsubhead
The results in this paper hold with similar proofs also for extended affine Weyl groups; to
simplify notation we do not treat this slightly more general case.

\head Contents\endhead
1. Weighted affine Weyl groups and the ring $J_*$.

2. Use of modular representations.

3. Folding.

4. A geometric interpretation of $P_{y,w;L}$.

\head 1. Weighted affine Weyl groups and the ring $J_*$\endhead
\subhead 1.1\endsubhead
Let $W$ be an irreducible
affine Weyl group with a given set $S$ of simple reflections
 assumed to have at least $2$ elements. Let $Q$ be the set of all translations in $W$, that is
 the set of all 
$t\in W$ such that the $W$-conjugacy class of $t$ is finite. It is known that
$Q$ is a normal free abelian subgroup of finite index of $W$.
We write the group operation in $Q$ as $+$.
We fix $s_0\in S$ such that $W$ is generated by $Q$ and by the finite
subgroup $W_0$ generated by $S_0=S-\{s_0\}$. (Such $s_0$ is said to be
``special''.)
Let $w\m|w|$ be the length function of $W$. Let
$Q^+=\{x\in Q;|sx|=|x|+1\text{ for any }s\in S_0\}$.
We have $W=\sqc_{x\in Q^+}W_0xW_0$. For any $x\in Q^+$ we denote by
$M_x$ the unique element in $W_0xW_0$ such that $|M_x|$ is maximal, or
equivalently, such that $|sM_x|=|M_x|-1=|M_xs|$ for all $s\in S_0$.
In particular, $M_0$ is the longest element in $W_0$. 
Let $L:W\to\NN$ be a weight function, that is a function such that
  $L(ww')=L(w)+L(w')$ whenever $w,w'$ in $W$ satisfy $|ww'|=|w|+|w'|$; we
assume that $L(s)>0$ for any $s\in S$.
Let $v$ be an indeterminate. Let $H$ be the $\QQ(v)$-vector space with basis
$\{T_w;w\in W\}$. We can regard $H$ as an associative
algebra in which $T_wT_{w'}=T_{ww'}$ if $w,w'$ in $W$ satisfy
$|ww'|=|w|+|w'|$ and $(T_s+v^{-L(s)})(T_s-v^{L(s)})=0$ for $s\in S$.
Let $\{c_w;w\in W\}$ be the basis of $H$ defined in \cite{L83a} and
\cite{L03, 5.2}. (See \cite{KL79} for the case $L=||$.)
We have $c_w=\sum_{y\in W}v^{-L(w)+L(y)}P_{y,w;L}T_y$ where
$P_{y,w;L}\in\ZZ[v^2]$ is zero for all but finitely many $y$
(see \cite{L03, 5.4}).
Let $\ca=\ZZ[v,v\i]$ and let $H_\ca$ be the $\ca$-submodule of
$H$ spanned by $\{T_w;w\in W\}$
or equivalently by $\{c_w;w\in W\}$. This is a subring of $H$.

We set $\p_L=\sum_{e\in W_0}v^{2L(e)}$.
For $x\in Q^+$ we set $\boc_x=\fra{v^{L(M_0)}}{\p_L}c_{M_x}\in H$.

\subhead 1.2\endsubhead
From \cite{L03, 8.6} we see that for $w\in W,x\in Q^+,y\in Q^+$ we have

$c_wc_{M_y}\in\sum_{u\in W;|us|=|u|-1\qua\frl s\in S_0}\ca c_u$,

$c_{M_x}c_w\in\sum_{u\in W;|su|=|u|-1\qua\frl s\in S_0}\ca c_u$.
\nl
It follows that 

$c_{M_x}c_{M_y}\in\sum_{u\in W;|su|=|u|-1=|us|\qua\frl s\in S_0}\ca c_u$
\nl
so that

(a) $c_{M_x}c_{M_y}=\sum_{z\in Q^+}\ti r_{x,y,z;L}c_{M_z}$
\nl
where $\ti r_{x,y,z;L}\in\ca$ is zero for all but finitely many $z$.
Hence

(b) $\boc_x\boc_y=\sum_{z\in Q^+}r_{x,y,z;L}\boc_z$
\nl
  where $r_{x,y,z;L}=v^{L(M_0)}\p_L\i\ti r_{x,y,z;L}$.

\subhead 1.3\endsubhead
Let $w\in W$. We write $c_w=\sum_{u\in W}l_uT_u$ with $l_u\in\ca$.
From \cite{L03, 6.6a} we see by induction on $|w|$ that:

(a) If $s\in S$, $|sw|=|w|-1$, then $c_w\in(T_s+v^{-L(s)})H_\ca$; in other
words, $l_{su}=v^{L(s)}l_u$ for any $u\in W$ such that $|su|=|u|+1$.
\nl
Assume now that $|w|=|M_0w|+|M_0|$. We show:

(b) for any $u\in W$ such that $|M_0u|=|u|+|M_0|$ and any $e\in W_0$ we have
$l_{eu}=v^{L(e)}l_u$.
\nl
We argue by induction on $|e|$. If $|e|=0$, there is nothing to prove.
Assume now that $|e|>0$. We have $e=se'$ for some $s\in S_0,e'\in W_0$ with
$|e|=|e'|+1$. By the induction hypothesis we have $l_{e'u}=v^{L(e')}l_u$.
We have $|se'u|=|e'u|+1$ (both sides are equal to $|se'|+|u|$). Using (a)
we have $l_{se'u}=v^{L(s)}l_{e'u}$ hence

$l_{eu}=l_{se'u}=v^{L(s)}v^{L(e')}l_u=v^{L(s)+L(e')}l_u=v^{L(se')}l_u=v^{L(e)}l_u$.
\nl
This proves (b).

In the setup of (b) we have
$$\align&c_w=\sum_{u\in W; |M_0u|=|M_0|+|u|}l_u\sum_{e\in W_0}v^{L(e)}T_eT_u\\&
=\sum_{e\in W_0}v^{L(e)}T_e\sum_{u\in W; |M_0u|=|M_0|+|u|}l_uT_u=
v^{L(M_0)}c_{M_0}\sum_{u\in W; |M_0u|=|M_0|+|u|}l_uT_u.\endalign$$

It follows that

(c) If $|w|=|M_0w|+|M_0|$, then $c_w\in c_{M_0}H_\ca$.
\nl
Similarly,

(d) if $w'\in W$ satisfies $|w'|=|w'M_0|+|M_0|$ then $c_{w'}\in H_\ca c_{M_0}$.
\nl
Taking $w'=M_x$, $w=M_y$ with $x,y $ in $Q^+$ we see from (c),(d) that

$c_{M_x}c_{M_y}\in H_\ca c_{M_0}c_{M_0}H_\ca\sub\p_LH_\ca$.
\nl
(We use that $c_{M_0}c_{M_0}\in\p_LH_\ca$.)
Combining this with 1.2(a) we see that $\ti r_{x,y,z;L}$ in 1.2 is in $\p_L\ca$.

If $x\in Q^+$ then from (c),(d) we see that

$c_{M_0}c_{M_x}=c_{M_x}c_{M_0}=v^{-L(M_0)}\p_Lc_{M_x}$.
\nl
(We use that $c_{M_0}c_{M_0}=v^{-L(M_0)}\p_Lc_{M_0}$.)

\subhead 1.4\endsubhead
From \cite{L03, 13.4}, for any $w,w'$ in $W$ we have

$T_wT_{w'}\in v^{L(M_0)}\sum_{w''\in W}\ZZ[v\i]T_{w''}$.
\nl
(In the case where $L=||$ this is proved in \cite{L85,\S7}; the proof for general
$L$ is entirely similar.) As in  \cite{L03, 13.5}, we deduce that for any $w,w'$ in $W$ we have

(a) $c_wc_{w'}=\sum_{w''\in W}h_{w,w',w''}c_{w''}$
\nl
(finite sum) where $h_{w,w',w''}=N_{w,w',w'';L}v^{L(M_0)}\mod v^{L(M_0)-1}$
with $N_{w,w',w'';L}\in\ZZ$.

Let $J$ be the free abelian group with basis $\{\t_w;w\in W\}$. We define a bilinear
multiplication $J\T J@>>>J$ by

$\t_w\t_{w'}=\sum_{w''\in W}N_{w,w',w'';L}\t_{w''}$,
\nl
(this is a finite sum.) It is known \cite{L03, 18.3} that this multiplication is associative
if the conditions in \cite{L03, 18.1} are satisfied.

Let $J_*$ be the subgroup of $J$ with basis $\{\t_{M_x};x\in Q^+\}$. From 1.2(b) we see that $J_*$
is closed under the multiplication in $J$; thus for $x,y$ in $Q^+$ we have

$\t_{M_x}\t_{M_y}=\sum_{z\in Q^+}N_{M_x,M_y,M_z;L}\t_{M_z}$,
\nl
(this is a finite sum.)

\proclaim{Theorem 1.5}(a) For $x,y$ in $Q^+$ we have
$$\boc_x\boc_y=\sum_{z\in Q^+}N_{M_x,M_y,M_z;L}\boc_z$$

(b) The subgroup $R$ of $H$ with $\ZZ$-basis $\{\boc_x;x\in Q^+\}$
is closed under multiplication in $H$.

(c) The isomorphism of abelian groups $R@>\si>>J_*$ given by $\boc_x\m\t_{M_x}$
is compatible with the multiplication. In particular, the multiplication in $J_*$
is associative.
\endproclaim
For $x,y,z$ in $Q^+$ we have
$$\align&r_{x,y,z;L}\p_Lv^{-L(M_0)}=\ti r_{x,y,z;L}=h_{M_x,M_y,M_z}=v^{L(M_0)}X\\&=
(\sum_{e\in W_0}v^{2L(e)})Y=(\sum_{e\in W_0}v^{-2L(e)})Y'\endalign$$
where $X\in\ZZ[v\i]$, $Y\in\ca$, $Y'\in\ca$. It follows that

$(\sum_{e\in W_0}v^{-2L(e)})\i X\in\ZZ[v,v\i]$.
\nl
Since $X\in\ZZ[v\i]$ and $\sum_{e\in W_0}v^{-2L(e)}\in1+v\i\ZZ[v\i]$, we have

$(\sum_{e\in W_0}v^{-2L(e)})\i X\in\ZZ[[v\i]]$;
\nl
but this is also in $\ZZ[v,v\i]$ hence it must be in $\ZZ[v\i]$. Thus

$(\sum_{e\in W_0}v^{-2L(e)})\i X\in\ZZ[v\i]$ that is,

(d) $r_{x,y,z;L}\in\ZZ[v\i]$.
\nl
From the definition of $c_w$, we have

$\bar h_{w,w',w''}=h_{w,w',w''}$
\nl
for any $w,w',w''$ in $W$, where $\bar{}:\ca@>>>\ca$ is the ring involution which
takes $v^n$ to $v^{-n}$ for any $n$. Using this and the fact that $\p_Lv^{-L(M_0)}$ is fixed by
$\bar{}$,
we see that the left hand side of (d) is fixed by $\bar{}$ hence is necessarily in $\ZZ$.

Taking the coefficient of $v^{L(M_0)}$ in the two sides of the equality

$r_{x,y,z;L}\p_Lv^{-L(M_0)}=h_{M_x,M_y,M_z}$
\nl
in which $r_{x,y,z;L}\in\ZZ$, we see that $r_{x,y,z;L}=N_{M_x,M_y,M_z}$. 
This completes the proof of (a). Now (b),(c) are immediate consequences of (a).

\subhead 1.6\endsubhead
The ring $R$ has unit element $\boc_0$ and is known to be commutative; it follows that the ring
$J_*$ has unit element $\t_{M_0}$ and is commutative.
In the case where $L=||$, 1.5(b) recovers a result in \cite{L83}. For general $L$, 1.5(b)
recovers a result in \cite{K05}. But the present proof is simpler than that in these references.

\subhead 1.7\endsubhead
In this subsection we assume that $L=||$. In this case the ring $J$ in \cite{L03, 18.3}
is associative. In \cite{L97} we have categorified $J$ to a monoidal tensor category
with simple objects indexed by $W$. In particular $J_*$ is categorified 
to a monoidal tensor category $\un J_*$. It is known that (as a consequence of 1.5(b))
$R$ can be also categorified to a monoidal category $\SS$ known as the ``Satake category''.
The ring isomorphism $R@>\si>>J_*$ in 1.5(c) gives rise to an
equivalence of monoidal categories $\SS@>\si>>\un J_*$.

\head 2. Use of modular representations\endhead
\subhead 2.1\endsubhead
In this section we assume that $L=||:W@>>>\NN$.
Let $\kk$ be an algebraically closed field of characteristic $p\ge0$.
Let $G$ be an adjoint semisimple group over $\kk$ with a fixed
pinning (involving a maximal torus $T$). We assume that the Weyl group of
$G$ is $W_0$, the lattice of roots of $G$ with respect to $T$ is
$Q$ and that $W=W_0Q$ is the affine Weyl group associated in the usual way
to the dual group $G^*$.
Then $Q^+$ is the set of dominant weights of $G$. For $x\in Q^+$ let $V_x$ be a
Weyl module of $G$ over $\kk$ with highest weight $x$; let $L_x$ be a simple
rational $G$-module withhighest weight $x$.

Let $\r\in Q_\RR=\RR\ot Q$ be half the sum of all positive roots of $G$.

Let $\ch$ be the set of hyperplanes
$H_{\cha,m}=\{x\in Q_\RR;\cha(x+\r)=mp\}$ for various
coroots $\cha:Q_\RR@>>>\RR$ and various $m\in\ZZ$.
(When $p=0$, $\ch$ consists of the hyperplanes $H_{\cha,0}$.)

\subhead 2.2\endsubhead
We now assume that $p$ is a prime number, $p\gg0$.
Following Verma \cite{Ve} we identify $W$ with the subgroup $W_p$
of the group of affine transformations of $Q_\RR$
generated by the reflections in the hyperplanes in $\ch$ which
preserve the set $\ch$.

Let $x\in Q^+$ be such that $x\n\cup_{\cha,m}H_{\cha,m}$ and
$\cha_0(x)\le p(p-h+2)$ where $\cha_0$ is the highest coroot and $h$ is the
Coxeter number. It is known \cite{AJS, KL94, KT95} that, as virtual $T$-modules, we have 

(a) $L_x=\sum_{y\in Z_x}(-1)^{|w_yw_x|}\dim(\cv_{w_y,w_x})V_y$,
\nl
where $Z_x$ is the set of all $y\in Q^+$ in the same $W_p$-orbit as $x$;
$w_x,w_y$ are certain well defined explicit elements of $W_p$; $\cv_{w_y,w_x}$
is a $\CC$-vector space of dimension $P_{w_y,w_x;||}(1)$ 
defined in terms of the stalks of the intersection cohomology complex
of an affine Schubert variety associated to $G^*$.

As shown in \cite{L17, comments to [53]}, from (a) with $x$ of the form $x=px',x'\in Q^+$
one can deduce that for $y'\in Q^+$ we have

(b) $P_{M_{y'},M_{x'};||}(1)=\dim(V_{x'}^{y'})$
\nl
where $V_{x'}^{y'})$ is the $y'$-weight space of $V_{x'}$.
(Note that in our case we have automatically $x\n\cup_{\cha,m}H_{\cha,m}$.)
This provides a new proof of one of the main results in \cite{L83}.

\subhead 2.3\endsubhead
In this subsection we assume that $p=0$.
Let $\AA$ be the subring of $\QQ(v)$ consisting of elements which
have no pole for $v=1$. Let $H_\AA$ be the $\AA$-submodule of
$H$ spanned by $\{T_w;w\in W\}$
or equivalently by $\{c_w;w\in W\}$. This is a subring of $H$.
We define a group homomorphism $\x$ from $H_\AA$ to the group ring $\QQ[W]$
by $\sum_wf_wT_w\m\sum_wf_w(1)w$; here $f_w\in\AA$. This is a ring
homomorphism. Recall that for $x\in Q^+$,
$P_{w,M_x;||}(1)$ depends only on the $(W_0,W_0)$ double coset of $w\in W$.
Hence
$$\align&\x(\boc_x)=\sha(W_0)\i\sum_{w\in W}P_{w,M_x;||}(1)w=
\sha(W_0)\i\sum_{x'\in Q^+}P_{M_{x'},M_x;||}(1)\sum_{w\in W_0x'W_0}w\\&=
\sha(W_0)\i\sum_{x'\in Q^+}\dim(V_x^{x'})\sum_{w\in W_0x'W_0}w.\endalign$$
(We have used 2.2(b).) We have also
$$\x(\boc_x)=\sha(W_0)\i\sum_{e\in Q}\dim(V_x^e)\sum_{a\in W_0}ae
           =\sha(W_0)\i\sum_{e'\in Q}\dim(V_x^{e'})\sum_{a\in W_0}e'a.\tag a$$
Indeed,
$$\align&\sum_{e\in Q}\dim(V_x^e)\sum_{a\in W_0}ae=\sum_{x'\in Q^+,(a,b)\in W_0\T W_0}
\fra{\dim(V_x^{x'})}{\sha(b'\in W_0;b'x'=x'b')} abx'b\i\\&=
\sum_{x'\in Q^+,w\in W_0x'W_0}\fra{\dim(V_x^{x'})\sha((b,c)\in W_0\T W_0,w=cx'b\i)}
{\sha(b'\in W_0;b'x'=x'b')} w\\&
=\sum_{x'\in Q^+,w\in W_0x'W_0}\fra{\dim(V_x^{x'})\sha((b,c)\in W_0\T W_0,cx'b\i=x')}
{\sha(b'\in W_0;b'x'=x'b')} w  \\&
  =\sum_{x'\in Q^+,w\in W_0x'W_0}\dim(V_x^{x'})w=\sha(W_0)\x(\boc_x)\endalign$$

and the first equality in (a) is established. The second equality in (a)
  follows the first by the substitution $e'=aea\i$.

\mpb
  
  Now let $x\in Q^+$, $y\in Q^+$. For $e''\in Q$ let $(V_x\ot V_y)^{e''}$
  be the $e''$-weight space of $V_x\ot V_y$. We have
$$\align&\x(\boc_x\boc_y)=\x(\boc_x)\x(\boc_y)\\&=
\sha(W_0)^{-2}\sum_{e\in Q}\dim(V_x^e)\sum_{a\in W_0}ae
\sum_{e'\in Q}\dim(V_y^{e'})\sum_{b\in W_0}e'b\\&
=\sha(W_0)^{-2}\sum_{(e,e')\in Q\T Q}\dim(V_x^e)\dim(V_y^{e'})\sum_{(a,b)\in W_0\T W_0}aee'b\\&
=\sha(W_0)^{-2}\sum_{e''\in Q}\dim(V_x\ot V_y)^{e''}\sum_{(a,b)\in W_0\T W_0}ae''b\\&
=\sha(W_0)^{-2}\sum_{e''\in Q,z\in Q^+}(V_z:V_x\ot V_y)\dim(V_z^{e''})
\sum_{(a,b)\in W_0\T W_0}ae''b\\&
=\sha(W_0)\i\sum_{e''\in Q,z\in Q^+}(V_z:V_x\ot V_y)\dim(V_z^{e''})
\sum_{a\in W_0}ae''.\tag b\endalign$$
Her $(V_z:V_x\ot V_y)$ is the multiplicity of $V_z$ in $V_x\ot V_y$.
On the other hand we have
$$\align&\x(\boc_x\boc_y)=\sum_{z\in Q^+}r_{x,y,z;||}\x(\boc_z)\\&
=\sha(W_0)\i\sum_{z\in Q^+,e''\in Q}r_{x,y,z;||}\dim(V_z^{e''})\sum_{a\in W_0}ae''.\endalign$$
Comparing with (b) we deduce
$$\sum_{z\in Q^+}(V_z:V_x\ot V_y)\dim(V_z^{e''})=\sum_{z\in Q^+}r_{x,y,z;||}\dim(V_z^{e''})$$
for any $e''\in Q$. Hence
$$\sum_{z\in Q^+}(V_z:V_x\ot V_y)V_z=\sum_{z\in Q^+}r_{x,y,z;||}V_z$$
in the Grothendieck group of representations of $G$. Since
  $(V_z)_{z\in Q^+}$ is a basis of this Grothendieck group, we see that

(c) $(V_z:V_x\ot V_y)=r_{x,y,z;||}$ 
\nl
for any $x,y,z$ in $Q^+$. Thus, we recover one of the main results in \cite{L83}.
    
\head 3. Folding\endhead
\subhead 3.1\endsubhead
In this section we assume that $W,S,s_0,W_0,Q,Q^+$ in 1.1
are such that $W$ is irreducible of simply laced type.
We asume given an automorphism $\s$ of $(W,S)$ of order
$\d\in\{2,3\}$ preserving $s_0$.
\nl
Let ${}'W=\{w\in W;\s(w)=w\}$.
For each $\s$-orbit $\co$ in $S$ let $s_\co$ be the longest element
in the subgroup of $W$ generated by the elements in $\co$.
Let ${}'S$ be the subset of ${}'W$ consisting of the elements $s_\co$
for various $\co$ as above. Note that $({}'W,{}'S)$ is an affine Weyl group.
Let $L:{}'W@>>>\NN$ be the restriction to ${}'W$ of the usual
length function of $W$; this is a weight function on
${}'W$.

We preserve the setup of 2.1. We assume that $G$ is simple of simply laced type.
We fix an automorphism of $G$ preserving the pinning of $G$ which induces the automorphism
$\s$ of $W$ considered above. This automorphism of $G$ is denoted again by $\s$.
If $x\in Q^+$ and $\s(x)=x$ then $\s:G@>>>G$ induces linear isomorphisms
$V_x@>>>V_x$,  $L_x@>>>L_x$ denoted again by $\s$ (they act as
identity on a highest weight vector). We have

$V_x=\op_{\th\in\kk^*_\d}V_{x,\th},L_x=\op_{\th\in\kk^*_\d}L_{x,\th}$  
\nl
where $\kk^*_\d=\{\th\in\kk^*;\th^\d=1\}$ and
$V_{x,\th},L_{x,\th}$ are the $\th$-eigenspaces of $\s$.

\subhead 3.2\endsubhead
We now assume that $p\gg0$ and that $x$ in 2.2(a) satisfies in addition $\s(x)=x$.
The proof of 2.2(a) is sufficiently functorial to imply that we have also
$$\sum_{\th\in\kk^*_\th}\ti\th L_{x,\th}=\sum_{y\in Z_x,\s(y)=y}
(-1)^{L(w_yw_x)}\tr(\s,\cv_{w_y,w_x})\sum_{\th\in\kk^*_\th}\ti\th V_{y,\th}\tag a$$
(equality in the representation ring of $T/\{\s(t)t\i;t\in T\}$
tensored with $\CC$; here $\th\m\ti\th$ is an imbedding of $\kk^*_\d$ into
$\CC^*$). Note that $\s(w_x)=w_x$ and that when $y\in Z_x$, $\s(y)=y$,
we have $\s(w_y)=w_y$, so that $\s$ acts naturally on $\cv_{w_y,w_x}$. We now substitute

(b) $\tr(\s,\cv_{w_y,w_x})=P_{w_y,w_x;L}(1)$
\nl
where $P_{w_y,w_x;L}$ is defined in terms of ${}'W$ and $L:{}'W@>>>\NN$ as in 3.1.
(See 4.5, 4.6.) We obtain

$$\sum_{\th\in\kk^*_\th}\ti\th L_{x,\th}=\sum_{y\in Z_x,\s(y)=y}
(-1)^{L(w_yw_x)}P_{w_y,w_x;L}(1)\sum_{\th\in\kk^*_\th}\ti\th V_{y,\th}\tag c$$
This is an extension of the character formula 2.2(a) to disconnected groups.
Note that the coefficients $P_{w_y,w_x;L}(1)$ are computable by an algorithm in
\cite{L03,\S6} (which is somewhat more involved than that for the unweighted case in \cite{KL79}).

\subhead 3.3\endsubhead
Note that $\s$ acts naturally on $G^*$.
Let ${}'G$ be the simply connected group over $\kk$ isogenous to
the dual group of the identity component of
the $\s$-fixed point set on $G^*$. By a theorem of Jantzen \cite{Ja}, the expression
$\sum_{\th\in\kk^*_\th}\ti\th V_{y,\th}$ in (c) can be expressed in terms of the character
of a Weyl module of ${}'G$. Using this one can deduce as in \S2 the analogues of
2.3(b), 2.4(c) with $(W,S,||)$ replaced by $({}'W,{}'S,L)$. (This recovers in our case a result in
\cite{K05}).

\subhead 3.4\endsubhead
Assume that $(W,S)$ is of (affine) type $A_2$ with $\s$ of order $2$. In this case $({}'W,{}'S)$
is of (affine) type $A_1$ and the values of $L|_{{}'S}$ are $1$ and $3$. In this case the ring $J_*$
associated to $({}'W,{}'S,L)$ in 1.5 is isomorphic together with its basis to the representation
ring of $SL_2(\CC)$ with its standard basis, see \cite{L03, 18.5}. This shows that the group
${}'G$ in 3.3 cannot be replaced by the corresponding adjoint group (even though
$G$ was adjoint). 

\subhead 3.5\endsubhead
In the setup of 3.1, 3.2  with $\kk=\CC$, we identify $W_0$ with the group $\cw_0$ of
affine transformations of $Q_\RR$ generated by the reflections in the (finitely many)
hyperplanes in $\ch$ and which preserve $\ch$. Let $\fg$ be the Lie algebra of $G$.
Let $x\in Q$ be such that $x\n\cup_{\cha}H_{\cha,0}$. Let $z\in Q$.
Then the Verma $\fg$-module $\VV_x$, its irreducible quotient $\LL_x$ 
and their $z$-weight spaces $\VV_x^z$, $\LL_x^z$ are defined.
It is known that the following equality (conjectured in \cite{KL79}) holds:

(a) $\dim\LL_x^z=\sum_{y\in\cz_x}(-1)^{|\o_y\o_x|}P_{\o_y,\o_x;||}(1)\dim\VV_y^z$,
\nl
where $\cz_x$ is the set of all $y\in Q$ in the same $\cw_0$-orbit as $x$;
$\o_x,\o_y$ are certain well defined explicit elements of $\cw_0$.

Now assume that $x,z$ are fixed by $\s$.
Then $\s:G@>>>G$ induces automorphisms of $\LL_x^z$ and of $\VV_x^z$ denoted again by $\s$. We have

(b) $\tr(\s,\LL_x^z)=\sum_{y\in\cz_x,\s(y)=y}(-1)^{L(\o_y\o_x)}
P_{\o_y,\o_x;L}(1)\tr(\s,\VV_y^z)$.
\nl
This follows from the proof of (a) in the same way as 3.2(c) follows from the proof of 2.2(a)
(using 4.5).

\head 4. A geometric interpretation of $P_{y,w;L}$\endhead
\subhead 4.1\endsubhead
Let $W_0$ be a (finite) Weyl group with a set $S_0$ of simple reflections and let
$\s:W_0@>>>W_0$ be an automorphism preserving $S_0$.
For each $\s$-orbit $\co$ in $S_0$ we denote by $\s_\co$ the longest element in the subgroup
of $W_0$ generated by the reflections in $\co$. 
Let ${}'W_0=\{w\in W_0;\s(w)=w\}$ and let ${}'S_0$ be the subset of ${}'W_0$ consisting of the
elements $s_\co$ for various $\co$ as above. Then ${}'W_0$ is a Weyl group with set
of simple reflections ${}'S_0$. Let $L:{}'W_0@>>>\NN$
be the restriction to ${}'W_0$ of the standard
length function of $W_0$; it is known that $L$ is a weight function on ${}'W_0$ so
that the Hecke algebra over $\ca$ 
with its bases $\{T_w;w\in{}'W_0\}$, $\{c_w;w\in{}'W_0\}$ can be defined as in 1.1
(in terms of ${}'W_0,{}'S_0,L$ instead of $W,S,L$). This Hecke algebra specialized at $v=\sqrt{q}$ with $q$ a prime power
is a $\CC$-algebra denoted by $H_{0,q;L}$.

For $w\in{}'W_0$ we write
$$c_w=\sum_{y\in{}'W_0}v^{-L(w)+L(y)}P_{y,w;L}T_y$$
where $P_{y,w;L}\in\ZZ[v^2]$.

For $w\in {}'W_0$ we have
$$T_wT_{w_0}=\sum_{y\in{}'W_0}v^{L(y)-L(w)}R_{y,w;L}T_{yw_0}\tag a$$
where $R_{y,w;L}\in\ZZ[v^2]$ is $0$ unless $y\le w$ and $w_0$ is
the longest element of $W_0$ (or ${}'W_0$).
Note the following inductive formulas for $R_{y,w;L}$. (Here $s\in S$.)
$$\align&R_{y,w;L}=R_{sy,sw;L}\text{ if }|sy|<|y|,|sw|<|w|;\\&
R_{y,w;L}=v^{2L(s)}R_{sy,w;L}+(v^{2L(s)}-1)R_{sy,sw}\text{ if }|sy|>|y|,|sw|<|w|.\tag b
\endalign$$
We have $P_{y,w;L}=0$ unless $y\le w$ and $P_{w,w;L}=1$. For $y,w$ in ${}'W_0$ we have
$$v^{2L(w)-2L(y)}\bP_{y,w;L}=\sum_{z\in{}'W_0}R_{y,z;L}P_{z,w;L}.\tag c$$

\subhead 4.2\endsubhead
Let $\kk$ be an algebraic closure of the finite prime field $\FF_p$.
Let $G$ be a simply connected semisimple group over $\kk$ with Weyl group $(W_0,S_0)$
and with a fixed pinning involving a maximal torus $T$ and a Borel subgroup $B$ containing $T$.
We fix an $\FF_p$-rational structure on $G$ (with Frobenius map $F:G@>>>G$) compatible with 
the pinning such that $T$ is split over $\FF_p$ hence $B$ is defined over $\FF_p$.
We consider an automorphism of $G$ preserving the pinning and compatible
with the $\FF_p$-structure;
it induces an automorphism of $W_0$ which we assume to be $\s$. This automorphism of $G$ is
denoted again by $\s$; we have $\s F=F\s$. Hence if $t\ge1$ then $F_t:=F^t\s=\s F^t$ is the
Frobenius map for a rational structure over the subfield $\FF_{p^t}$ with $p^t$ elements of $\kk$.
Let $\cb$ be the variety of Borel subgroups of $G$. Note that $F_t$ acts naturally on $\cb$
and defines a Frobenius map on $\cb$. 
We say that $B_1,B_2$ in $\cb$ are opposed if $B_1\cap B_2$ is a maximal torus.
We define $B^*\in\cb$ by the conditions that $B\cap B^*=T$.
For $B_1,B_2$ in $\cb$ let $pos(B_1,B_2)\in W_0$ be the relative position of $B_1,B_2$.
For $w\in W_0$ we set $\cb_w=\{B'\in\cb;pos(B,B')=w\}$.
For $y\in W_0$ we define ${}^yB\in\cb$ by the conditions $T\sub {}^yB, {}^yB\in\cb_y$;
we define ${}^yB^*\in\cb$ by the conditions $T\sub {}^yB^*, pos(B^*,{}^yB^*)=y$.
Let $\bar\cb_w$ be the closure of $\cb_w$ in $\cb$.
For $y\in W_0$ we set $A^y=\{B'\in\cb;B',{}^yB^*\text{ opposed}\}$.

For any algebraic variety $X$ of pure dimension let $IC(X)$ be the intersection cohomology complex
of $X$ with coefficients in $\bbq$ (with $l$ a prime $\ne p$). Let $\HH^i(X)$ (resp. $\HH^i_c(X)$) be the $i$-th
cohomology (resp. $i$-th cohomology with compact support) of $X$ with coefficients in $IC(X)$.
For $x\in X$ let $\ch^i_x(X)$ be the stalk at $x$ of the $i$-th cohomology sheaf of $IC(X)$.

The following result gives a geometric interpretation of $P_{y,w;L}$ (stated without proof
in \cite{L83a, (8.1)}) extending the already known case where $\s=1$ considered in \cite{KL80};
see also \cite{L03, \S16}.

\proclaim{Theorem 4.3}Let $y\in{}'W_0,w\in{}'W_0$ be such that $y\le w$. We have
$$P_{y,w;L}=\sum_{i\text{ even}} \tr(\s,\ch^i_{{}^yB}(\bar\cb_w)v^{2i}.$$
\endproclaim
(Note that $\s$ acts naturally on $\ch^i_{{}^yB}(\bar\cb_w)$.)
The proof will use the following result (analogous to \cite{KL79, A4(a)}).

\proclaim{Lemma 4.4} Let $y\in{}'W_0,z\in{}'W_0$ be such that $y\le z$. We have
$\sha((\cb_z\cap A^y)^{F_t})=R_{y,z;L}(p^t)p^{tL(y)}$.
\endproclaim
Let $\cf$ be the vector space of functions $\cb^{F_t}@>>>\CC$. Then $\cf$ is an $H_{0,p^t;L}$-module
in which for $w\in{}'W_0$ and $f\in\cf$ we have $T_wf=f'$ where for $B'\in\cb^{F_t}$ we have
$f'(B')=p^{-tL(w)/2}\sum_{B''\in\cb^{F_t};\pos(B',B'')=w}f(B'')$.
Applying the equality 4.1(a) to $f\in\cf$ and evaluating at $B$ we see that for $w\in{}'W_0$ we have
$$\align&\sum_{y'\in{}'W_0}p^{t(L(y')-L(z))/2}R_{y',z;L}(p^t)p^{t(L(y')-L(w_0))/2}  
\sum_{C\in\cb^{F_t};pos(B,C)=y'w_0}f(C)\\&
=p^{-tL(z)/2}\sum_{B''\in\cb^{F_t};pos(B,B'')=z}p^{-tL(w_0)/2}\sum_{C\in\cb^{F_t};pos(B'',C)=w_0}f(C).
\tag a\endalign$$
We now take $f$ to be the function equal to $1$ at $C_0={}^{yw_0}B$ and equal to $0$ on
$\cb^{F_t}-\{C_0\}$. We obtain
$$\align&\sha(B''\in\cb^{F_t};pos(B,B'')=z,pos(B'',C_0)=w_0)p^{-tL(z)/2}\\&=
p^{t(L(y)-L(z)+L(y))/2}R_{y,z;L}(p^t),\tag b\endalign$$
that is
$$\sha(\cb_z\cap A^y)^{F_t}=R_{y,z;L}(p^t)p^{tL(y)}.$$
The lemma is proved.

\subhead 4.5\endsubhead
We now prove the theorem. When $y=w$ the result is obvious.
We can assume that $y<w$ and that

(a) the result is true when $y,w$ is replaced by $z,w$ with $z\in{}'W_0$
such that $y<z\le w$.
\nl
Applying the Grothendieck-Lefschetz fixed point formula for $F_t$
on the $F_t$-stable open subvariety $\bar\cb_w\cap A^y$ of $\bar\cb_w$ we obtain
$$\tr(F_t,\sum_i(-1)^i\HH^i_c(\bar\cb_w\cap A^y))=\sum_{z\in W_0;y\le z\le w}
\sum_{B'\in(\cb_z\cap A^y)^{F_t}}\tr(F_t,\sum_i(-1)^1\ch^i_{B'}(\bar\cb_w)).$$
Now the fixed point set $(\cb_z\cap A^y)^{F_t}$ is empty unless $\s(z)=z$. For such $z$ we
apply Lemma 4.4 and we obtain
$$\align&\tr(F_t,\sum_i(-1)^i\HH^i_c(\bar\cb_w\cap A^y))\\&=\sum_{z\in{}'W_0;
y\le z\le w}
R_{y,z;L}(p^t)p^{tL(y)}\tr(F_t,\sum_i(-1)^i\ch^i_{{}^zB}(\bar\cb_w)).\endalign$$
By Poincar\'e duality on $\bar\cb_w\cap A^y$ we have
$$\tr(F_t,\sum_i(-1)^i\HH^i_c(\bar\cb_w\cap A^y))
=p^{tL(w)}\tr(F_t\i,\sum_i(-1)^i\HH^i(\bar\cb_w\cap A^y)).$$
Using \cite{KL80, 4.5(a), 1.5} we have
$$\tr(F_t\i,\sum_i(-1)^i\HH^i(\bar\cb_w\cap A^y))=\tr(F_t\i,\sum_i(-1)^i\ch^i_{{}^yB}(\bar\cb_w)),$$
so that 
$$\align&p^{tL(w)}\tr(F_t\i,\sum_i(-1)^i\ch^i_{{}^yB}(\bar\cb_w))\\&
=\sum_{z\in {}'W_0;y\le z\le w}
R_{y,z;L}(p^t)p^{tL(y)}\tr(F_t,\sum_i(-1)^i\ch^i_{{}^zB}(\bar\cb_w)).\endalign$$
By \cite{KL80, 4.2} we have $\ch^i_{{}^yB}(\bar\cb_w)=0$ if $i$ is odd while if $i$ is even
the eigenvalues of $F^t$ on $\ch^i_{{}^yB}(\bar\cb_w)$ are equal to $p^{it/2}$. It follows that
$$\align&p^{tL(w)}\sum_{i\text{ even}} p^{-it/2}\tr(\s\i,\ch^i_{{}^yB}(\bar\cb_w))\\&=
\sum_{z\in {}'W_0;y\le z\le w}R_{y,z;L}(p^t)p^{tL(y)}
\sum_{i\text{ even}}p^{it/2}\tr(\s,\ch^i_{{}^zB}(\bar\cb_w)).\endalign$$
Since this holds for $t=1,2,\do$ we can replace $p^t$ by $v^2$ where $v$ is an indeterminate
and we get an equality in $\bbq[v,v\i]$:
$$\align&v^{2L(w)}\sum_{i\text{ even}} v^{-i}\tr(\s\i,\ch^i_{{}^yB}(\bar\cb_w))\\&=
\sum_{z\in{}'W_0;y\le z\le w}R_{y,z;L}v^{2L(y)}
\sum_{i\text{ even}}v^i\tr(\s,\ch^i_{{}^zB}(\bar\cb_w)).\endalign$$
Using the induction hypothesis (a) we obtain
$$\align&v^{2L(w)}\sum_{i\text{ even}} v^{-i}\tr(\s\i,\ch^i_{{}^yB}(\bar\cb_w))-
v^{2L(y)}\sum_{i\text{ even}}v^i\tr(\s,\ch^i_{{}^yB}(\bar\cb_w))\\&=
\sum_{z\in {}'W_0;y<z\le w}R_{y,z;L}v^{2L(y)}P_{z,w;L}.\endalign$$
Using 4.1(c), the right hand side of this equality is
$$v^{2L(w)}\bP_{y,w;L}-v^{2L(y)}P_{y,w;L}.$$
Thus we have 
$$\align&
v^{L(w)-L(y)}(\sum_{i\text{ even}} v^{-i}\tr(\s\i,\ch^i_{{}^yB}(\bar\cb_w))-\bP_{y,w;L})\\&=
v^{L(y)-L(w)}(\sum_{i\text{ even}}v^i\tr(\s,\ch^i_{{}^yB}(\bar\cb_w))-P_{y,w;L}).\tag b\endalign$$
By the known properties of $\ch^i$, in both sides of (b) we can assume that
$i<\dim\cb_w-\dim\cb_y=L(w)-L(y)$. Moreover, we have $v^{L(w)-L(y)}\bP_{y,w;L}\in v\ZZ[v]$
and $v^{L(y)-L(w)}P_{y,w;L}\in v\i\ZZ[v\i]$.
Thus the left hand side of (b) is in $v\bbq[v]$ while the
right hand side of (b) is in $v\i\bbq[v\i]$.
We see that both sides of (b) are zero. The theorem is proved.

\subhead 4.6\endsubhead
The proof in 4.5 is written in such a way that it remains valid in the affine case
so that it gives an analogous geometric interpretation for $P_{y,w;L}$ with $y,w$ in ${}'W$ where
${}'W,L$ are as in 3.1. In this case $B$ is an Iwahori subgroup and $B^*$ is an anti-Iwahori
subgroup (opposed to $B$) as in \cite{KL80,\S5}.
The definition of $A^y$ still makes sense; it is the set of Iwahori subgroups opposed to a certaain
fixed anti-Iwahori subgroup. Now $R_{y,w;L}$ as defined by 4.1(a) does not make sense in the affinne
case; instead one can use the inductive definition in 4.1(b). With this definition, the
analogue of 4.1(c) remains valid; Lemma 4.4
remains valid but it is now proved by an (easy)  induction on $|z|$.

\subhead 4.7\endsubhead
{\it Erratum to \cite{L83}.} On p.212, line -3, the definition of $K^1$ should be

$K^1=\{x\in(1/|W|)\ZZ[\tW_a];wx=x,xw=x\qua\frl w\in W\}$.
\nl
On p.212, line -1, the definition of $J^1$ should be

$J^1=\{x\in\ZZ[\tW_a];wx=(-1)^{l(w)}x,xw=x\qua\frl w\in W\}$.
\nl
On p.212, line -9, the definition of $K$ should be 

$K=\{x\in(1/\cp)H;T_wx=q^{l(w)}x,xT_w=q^{l(w)}x\qua\frl w\in W\}$.
\nl
On p.212, line -7, the definition of $J$ should be 

$J=\{x\in H;T_wx=(-1)^{l(w)}x,xT_w=q^{l(w)}x\qua\frl w\in W\}$.


\widestnumber \key{KL94}
\Refs
\ref\key{AJS}\by H.H.Andersen, J.C.Jantzen and W.Soergel\paper Representations of quantum groups at
a $p$-th root of unity and of semisimple groups in characteristic $p$: independence of $p$\jour
Ast\'erisque\vol220\yr1994\pages1-321\endref
\ref\key{Ja}\by J.C.Jantzen\paper Darstellungen halbeinfacher algebraischer
Gruppen\jour Bonn.Math.Sch., no.67\yr1973\endref
\ref\key{KT}\by M.Kashiwara and T.Tanisaki\paper Kazhdan-Lusztig conjecture for affine Lie
algebras with negative level\jour Duke J. Math.\vol77\yr1995\pages21-62\endref
\ref\key{KL79}\by D.Kazhdan and G.Lusztig\paper Representations of Coxeter groups and Hecke algebras
\jour Inv. Math.\vol53\yr1979\pages 165-184\endref
\ref\key{KL80}\by D.Kazhdan and G.Lusztig\paper Schubert varieties and Poincar\'e duality\inbook Proc.
Symp. Pure Math. 36\publ Amer. Math. Soc.\yr1980\pages185-203\endref
\ref\key{KL94}\by D.Kazhdan and G.Lusztig\paper Tensor structures arising from affine Lie algebras, IV
\jour J. Amer. Math. Soc.\vol7\yr1994\pages383-453\endref
\ref\key{K05}\by F.Knop\paper On the Kazhdan-Lusztig basis of a spherical Hecke algebra\jour
Represent.Th.\vol9\yr2005\pages417-425\endref
\ref\key{L83}\by G.Lusztig\paper Singularities, character formulas and a $q$-analog of weight
multiplicities\jour Ast\'erisque\vol101-102\yr1983\pages 208-229\endref
\ref\key{L83a}\by G.Lusztig\paper Left cells in Weyl groups\inbook Lie groups representations\bookinfo
LNM 1024\publ Springer Verlag\yr1983\pages99-111\endref
\ref\key{L85}\by G.Lusztig\paper Cells in affine Weyl groups\inbook Algebraic groups and related
topics \bookinfo Adv. Stud. Pure Math. 6\publ North-Holland and Kinokuniya\yr1985\pages 255-287
\endref
\ref\key{L97}\by G.Lusztig\paper Cells in affine Weyl groups and tensor categories\jour
Adv. Math.\vol129\yr1997\pages 85-98\endref
\ref\key{L03}\by G.Lusztig\book Hecke algebras with unequal parameters\bookinfo CRM Monograph
Ser.18\publ Amer. Math. Soc.\yr2003\endref 
\ref\key{L17}\by G.Lusztig\paper Comments on my papers\jour arxiv:1707.09368\endref
\ref\key{Ve}\by D.N.Verma\paper Role of affine Weyl groups in the representation theory of algebraic
Chevalley groups and their Lie algebras\inbook Lie Groups and representations, Budapest\yr1975\endref
\endRefs
\enddocument